\documentclass{amsart}
\usepackage{amsmath}
\usepackage{framed,comment,enumerate}
\usepackage[pdftex]{graphicx}
\usepackage{epstopdf}
\usepackage{xcolor, framed}
\usepackage{color}
\def\R {\mathbb{R}}
\def\eps{\varepsilon}
\def\OP{obstacle problem }
\def\FB{free boundary }

\def\Qua{\mathcal{Q}}

\def\Sph{\mathbb{S}^{d-1}}

\def\PosS{\{u>0\}}
\def\Singu{\Sigma(u)}

\newtheorem{prop}{Proposition}[section]
\newtheorem{thm}{Theorem}[section]

\newtheorem{lem}{Lemma}[section]
\theoremstyle{definition}

\numberwithin{equation}{section}

\title{A brief survey on the obstacle problem }

\author{Hui Yu}
\address{Department of Mathematics,	Columbia University, New York, USA}
\email{ huiyu@math.columbia.edu}
\begin{document}

\begin{abstract}
We discuss some regularity issues in the study of the obstacle problem.  In particular, we present a recent result by O.~Savin and the author on the regularity of the singular set for the obstacle problem with a fully nonlinear elliptic operator.  

This survey is based on a lecture by the author at the 8th International Congress of Chinese Mathematicians.
\end{abstract}

\maketitle

\section{Introduction}
Free boundary problems arise in the study of physical systems involving several distinct phases. These lead to equations with discontinuities along interfaces between the phases. In contrast with equations with \textit{prescribed} discontinuities, the key feature in free boundary problems is that the locations of these discontinuities are part of the \textit{unknown}. One of the goals is to understand the regularity of these unknown interfaces, the so-called free boundaries. 

The past few decades witnessed developments in many free boundary problems. Few problems, however, have achieved the status of the obstacle problem. Apart from its numerous direct applications, ideas and techniques originally developed for the obstacle problem have been adapted to many other problems. In this sense,  the obstacle problem is arguably the  archetypical free boundary problem. 

In this brief survey, we focus on  the obstacle problem.  We begin with the classical obstacle problem, that is, the obstacle problem with the Laplacian operator.  Although this problem has a long history, there have been exciting developments in the past few years, especially on the regularity of the singular part of the free boundary.  Then we move on to the obstacle problem with a fully nonlinear elliptic operator.  In particular, we discuss the recent resolution of the regularity of the singular set in a joint work of O.~ Savin and the author \cite{SY}.

This is a very brief survey. Instead of trying to be comprehensive, the intention is to highlight only a few results that appear important in the development of the general theory. Consequently, many interesting contributions are omitted. Fortunately, there are many surveys and books that paint a much more complete picture of the subject, for instance, Figalli \cite{F}, Petrosyan-Shahgholian-Uraltseva \cite{PSU}, and Ros-Oton \cite{R}.

\section{The classical obstacle problem}
The classical obstacle problem models the height of an elastic membrane being pushed towards an impenetrable obstacle.  Suppose that $\Omega$ is a domain in $\R^d.$ Along its boundary, an elastic membrane is fixed at height $x_{d+1}=g$ for some function $g:\partial\Omega\to[0,+\infty).$ Inside the domain we push the membrane downward (toward the hyperplane $\{x_{d+1}=0\}$) with constant force $1$.  At the level of the hyperplane $\{x_{d+1}=0\}$, there is an impenetrable obstacle. This forces the membrane to stay inside the region $\{x_{d+1}\ge 0\}.$

If we denote the height of the membrane by a function $u$, then this function solves the following equation, the so-called \textit{obstacle problem}\footnote{Although it is motivated by a very simple physical situation,  this equation also appears in numerous other problems, for instance, the melting of ice, fluid filtration, Hele-Shaw flows, and mathematical finance. For many of these applications, see the wonderful survey by Ros-Oton \cite{R} or the beautiful book by Petrosyan-Shahgholian-Uraltseva \cite{PSU}.}:\begin{equation}\label{ObstacleProblem}
\begin{cases}
\Delta u=\chi_{\PosS}&\\ u\ge 0&\end{cases}\text{ in $\Omega.$}
\end{equation}
Here we use the standard notation $\chi_E$ to denote the characteristic function of a set $E$.

The most interesting feature of this equation is the jump in the right-hand side. This arises from the different physical phases, depending on whether the membrane is above the obstacle or in contact with the obstacle. In the \textit{non-contact set} $\PosS$, the shape of the membrane is determined by the balance of its elasticity and the downward force, that is, $\Delta u=1.$  In the \textit{contact set} $\{u=0\},$ since the membrane cannot penetrate the obstacle, it has to be completely flat. 

This discontinuity occurs along the interface $\partial\PosS$, the \textit{free boundary} in this problem. Note that the location of this free boundary depends on the solution $u$ and is part of the unknown. 

The goal is to understand the regularity of the solution $u$ as well as the \FB $\partial\PosS$.  We focus on the interior regularity. 

\subsection{Regularity of the solution}
The first step is to understand the regularity of the solution to \eqref{ObstacleProblem}. 

Since the Laplacian of $u$ is  bounded, standard elliptic theory implies that the Hessian of $u$, $D^2u$, is in $\mathcal{L}^p_{loc}(\Omega)$ for all finite $p.$  An application of Sobolev embedding gives $C^{1,\alpha}_{loc}$ regularity of $u$ for all $\alpha\in(0,1).$ 

Since the Laplacian is not continuous, the solution is not in $C^2$. The \textit{optimal regularity} is $C^{1,1}_{loc}$, that is, the Hessian $D^2u$ is locally  bounded. This was established by Br\'ezis-Kinderlehrer \cite{BK}:
\begin{thm}\label{OptimalRegularity}
The solution $u$ to \eqref{ObstacleProblem} is in $C^{1,1}_{loc}(\Omega).$ 
\end{thm} 

For the regularity of the free boundary, it is crucial to understand the behavior of the solution near a point on the free boundary.  

To this end, we observe that in the interior of the contact set $\{u=0\}$, the Hessian satisfies $D^2u=0$. Inside the non-contact set $\PosS$, we have $\Delta u=1$, and that  $\Delta(D^2u)=0$. Since the function $u$ achieves its absolute minimum along the free boundary, in a weak sense we have $D^2u\ge0$ along $\partial\PosS.$ 

Combining these, it can be shown that the solution has precise quadratic behavior around a \FB point \cite{C2}:
\begin{prop}\label{QuadraticGrowth}Suppose that $u$ is a solution to \eqref{ObstacleProblem} with $0\in\partial\PosS$, then  $$cr^2\le \sup_{B_r}u\le Cr^2,$$ for all $r>0$ with $B_{2r}\subset\Omega.$ The constants $c$ and  $C$ depend only on the dimension $d.$
\end{prop} 

\subsection{Regular part in the  \FB}
The next step is to understand the regularity of the \FB $\partial\PosS$. 

Although there were results in two dimensions \cite{Sak1, Sak2, Sch}, these depend on complex-variable techniques and cannot be generalized to higher dimensions.  The first breakthrough in general dimensions was due to Caffarelli \cite{C1}. Here we follow the modern interpretation  \cite{C2}, where Caffarelli showed the following:
\begin{thm}\label{CaffarelliResult}
Let $u$ be a solution to the obstacle problem \eqref{ObstacleProblem}. 

The free boundary decomposes into two pieces, the regular part and the singular part, $$\partial\PosS\cap\Omega=Reg(u)\cup\Singu.$$

The regular part $Reg(u)$ is relatively open in $\partial\PosS$, and is locally an analytic hypersurface. 
\end{thm} 
To achieve this, Caffarelli introduced the technique of blow-up analysis to \FB problems. This has been the paradigm in the study of \FB regularity ever since.  

Roughly, there are two steps in this paradigm:
\begin{enumerate}
\item{\textit{Step 1: Rescale and blow up.} In this step, we study rescaled solutions and their limits. This allows us to magnify around a \FB point. Since we push the influence of boundary data to infinity,  we  obtain simpler objects as limits of rescalings, the so-called \textit{blow-up profiles}. }
\item{\textit{Step 2: Transfer information to the original solution.} At small scales, our solution resembles the blow-up profiles.  In this step, we quantify this resemblance and show that the \FB of our solution inherits regularity from the \FB of the blow-up profiles.}
\end{enumerate}

We illustrate these two steps in the obstacle problem. To simplify notations, we  assume that \textit{the origin is on the free boundary}, that is, $$0\in\partial\PosS.$$

Proposition \ref{QuadraticGrowth} says that the solution grows quadratically near a \FB point. Consequently, the quadratic rescaling is the only reasonable rescaling.  For $r>0$, we define the \textit{rescaled solution} \begin{equation*}\label{Rescaling}u_r(x)=\frac{1}{r^2}u(rx).\end{equation*} When the parameter $r$ is small, we are magnifying around the \FB point and pushing boundary data on $\partial\Omega$ to infinity. This simplifies the situation when $r\to0.$

The precise bounds in Proposition \ref{QuadraticGrowth} gives compactness and non-degeneracy for the family of rescalings.  In particular we have the following:
\begin{prop}\label{Compactness}
There is a sequence of $r_j\to 0$ such that for some function $u_0$ we have $$u_{r_j}\to u_0 \text{ locally uniformly in $C^{1,\alpha}(\R^d).$}$$
\end{prop}
This limit $u_0$ solves the obstacle problem \eqref{ObstacleProblem} in the entire $\R^d$ with $0\in\partial\PosS.$ Actually, we have a complete classification of the possible shapes of $u_0$:
\begin{prop}\label{Classification}
Let $u_0$ be as in Proposition \ref{Compactness}. 

Then $u_0$ is either a half-space solution of the form $$u_0(x)=\frac{1}{2}\max\{x\cdot e,0\}^2 \text{ for some $e\in\Sph$,}$$ or a parabola solution of the form $$u_0(x)=\frac{1}{2}x\cdot Ax$$ for some matrix $A$ with $A\ge 0$ and $trace(A)=1$.
\end{prop} 

The blow-up profile $u_0$ is much simpler than a general solution. Its free boundary is always flat. In the case of a half-space solution, the free boundary is the hyperplane perpendicular to the direction $e$. In the case of a parabola solution, the free boundary is the linear space $ker(A).$ 

This completes Step 1 in the blow-up analysis. 

In the next step, we transfer information from the blow-up profile $u_0$ back to our original solution $u$. In particular, we want to establish that if $u_0$ is a half-space solution, then the \FB of $u$ is similar to a hyperplane. 

However, the convergence in Proposition \ref{Compactness} only holds for a \textit{particular subsequence $r_j\to0,$} not the full sequence $r\to 0.$ It is not even clear if $u_{r}$ could converge to a half-space solution along some $r_j\to0$, but converge to a parabola solution along some other $r'_j\to0.$ 

This possible dependence on subsequences is one of the main difficulties in the study of \FB problems.

Note that for a half-space solution $\frac{1}{2}\max\{x\cdot e,0\}^2$, the contact set $\{u_0=0\}$ is a half space $\{x\cdot e\le 0\}$. For a parabola solution, the contact set is of codimension at least $1$, and has zero measure. This information can be transferred to the original solution. Around a regular point, there is true contact in a set of `full measure'. Around a singular point, there is only \textit{tangential contact} between the membrane and the obstacle.

This gives the following geometric characterization of the two possibilities in  Proposition \ref{Classification}:
\begin{prop}\label{GeometricCharacterization}
Let $u$ be  a solution to \eqref{ObstacleProblem} with $0\in\partial\PosS.$  

We have the following dichotomy:
\begin{enumerate}
\item{If $\limsup_{r\to 0}\frac{|B_r\cap\{u=0\}|}{r^d}>0$, then all blow-up profiles are half-space solutions; and }
\item{If $\limsup_{r\to 0}\frac{|B_r\cap\{u=0\}|}{r^d}=0$, then all blow-up profiles are parabola solutions.}\end{enumerate}
\end{prop} 
This characterization gives \textit{the uniqueness of type} of blow-up profiles.  

In particular, we can decompose the \FB depending on whether blow-ups are half-space solutions or parabola solutions. A point is called a \textit{regular point} if blow-ups are half-space solutions. Otherwise it is a \textit{singular point}.

The characterization in Proposition \ref{GeometricCharacterization} shows that at a regular point, the membrane contacts the obstacle `in full measure'. This is a very stable situation. It can be shown that if our solution $u$ is close to a half-space solution, then the free boundary $\partial\PosS$ is close to a hyperplane and is `almost flat'. An improvement of flatness argument gives analyticity of the regular part as in Theorem \ref{CaffarelliResult}.
\subsection{Instability of the singular part}The regularity of the regular part  depends crucially on the stability of the \FB of half-space solutions. Take for example $$u_0(x)=\frac{1}{2}\max\{x_1,0\}^2.$$ If we perturb this solution slightly on $\partial B_1$, then in $B_{1/2}$ the perturbed solution still contacts the obstacle `in full measure'. In particular, we still see a regular point. This observation allows us to show that if our solution $u$ is close to $u_0$, then the free boundary $\partial\PosS$ is close to $\partial\{x_1>0\}.$

This stability does not hold for parabola solutions. 

For instance, for the parabola solution $$u_0(x)=\frac{1}{2}x_1^2,$$ the contact set is the same as its free boundary $\{x_1=0\}.$ All points on this \FB are singular points.  

If we add a small positive constant $\eps$ to the data on $\partial B_1$, then  the perturbed solution $\frac{1}{2}x_1^2+\eps$ is strictly positive in the entire $B_1$. The free boundary disappears.  If we subtract $\eps$ from the data on $\partial B_1$, then the perturbed solution is of the form $\frac{1}{2}\max\{|x_1|-c\eps,0\}^2$. The contact set is a strip of positive width, and all the free boundary points become regular points.

Due to this instability, the study of the singular part is much more challenging than the regular part. 
\subsection{Regularity of the singular part in the free boundary}\label{SingularPartClassical}
Already in two dimensions, there is an example where the singular part is a Cantor set in the line $\{x_1=0\}$ \ \cite{Sch}.  As a result,  the singular part is not a smooth manifold in general. The best we can hope for is to show that \textit{the singular part is covered by manifolds with some regularity.}  

To simplify notations, let's define \textit{the space of parabola solutions}:
\begin{equation}
\Qua:=\{p=\frac{1}{2}x\cdot Ax| A\ge0, trace(A)=1\}.
\end{equation} 

In the following, we assume that \textit{the origin is a singular point on the free boundary,} that is, $$0\in\Singu.$$

By definition of the singular set, along some subsequence $r_j\to0$, the rescaled solutions, $u_{r_j}$, satisfy $$u_{r_j}\to p=\frac{1}{2}x\cdot Ax$$ for some $p\in\Qua.$ 

The natural strategy is to show that the singular set is covered by a manifold with tangent space $ker(A)$ at $0$. And in particular, around $0$ we should expect the manifold to have the same dimension as $dim(ker(A)).$ 

For this to work, however, it requires to show that the limit $p$ does not depend on the particular subsequence $r_j\to 0$. We need to show that $p$ is the \textit{unique} blow-up profile and $$u_r\to p$$ for the full $r\to0. $ This guarantees that the tangent space is unique, and that the manifold is differentiable.  If we can further compare  blow-ups at nearby points, then we  gain further regularity of this manifold. 

Due to the instability of the \FB for parabola solutions,  even the uniqueness of the blow-up is subtle.  The technical tools are various monotonicity formulae. These formulae allow quantification of the convergence of re-scaled solutions. They allow estimate of the form 
\begin{equation}\label{QuadraticApproximation}|u-p|\le\sigma(r)r^2 \text{ in $B_r$ for all $r>0$}\end{equation} for some modulus of continuity $\sigma. $  

This implies the uniqueness of the blow-up profile $p$. In particular, both the tangent space at $0$ and the dimension of the manifold are well-defined.  

The quantified modulus of continuity gives comparison of blow-up profiles at nearby points.  Suppose $x$ and $y$ are two singular points, and that $p_x$ and $p_y$ are the corresponding blow-up profiles. Then \eqref{QuadraticApproximation} implies $$|D^2p_x-D^2p_y|<\sigma(|x-y|).$$ Whitney's lemma implies that these parabolas at various points are the second order Taylor expansion of some function $f$. Moreover, the Hessian $D^2f$ has modulus of continuity $\sigma.$

Since $D^2f(0)=D^2p=A$, if $ker(A)$ is of dimension $k$, then $D^2f(0)$ contains a $(d-k)$-by-$(d-k)$ invertible submatrix, say $(\frac{\partial^2}{\partial x_i\partial x_j}f(0))_{1\le i,j\le d-k}$. Consequently, the map $\nabla' f=(\frac{\partial}{\partial x_1}f,\dots, \frac{\partial}{\partial x_{d-k}}f):\R^d\to\R^{d-k}$ is full-rank.  As a result, $(\nabla'f)^{-1}(0)$ is a manifold of dimension $k$ in a neighborhood of $0$.  This is the desired covering manifold of the singular set around the origin.

Implicit function theorem implies that this manifold has the same regularity as $\nabla 'f$. Since $D^2f$ has modulus of continuity $\sigma$, $\nabla'f$ is $C^1$. Its derivatives have modulus of continuity $\sigma.$ To get better regularity of the covering manifold, we need to get better control of the modulus of continuity $\sigma$ in \eqref{QuadraticApproximation}. 

To do this, Caffarelli's original treatment \cite{C2} depends on the Alt-Caffarelli-Friedman formula \cite{ACF}. Further developments by Weiss \cite{W} and Monneau \cite{M} rely on two monotonicity formulae now bearing their names. Recently, Colombo-Spolaor-Velichkov \cite{CSV} improved the results by further quantifying the analysis of Weiss monotonicity formula. The current best result on the singular set is due to Figalli-Serra \cite{FSe}. By introducing Almgren's monotonicity formula \cite{Alm} to the study of the singular set, they were able to show the following\footnote{The result by Figalli-Serra is more precise than what is stated here. For their complete result, consult \cite{FSe}.}
\begin{thm}\label{FSeResult}
Let $u$ be a solution to the obstacle problem \eqref{ObstacleProblem}.

The set of singular points stratifies $$\Singu=\cup_{k=0}^{d-1}\Sigma^k(u).$$

The highest stratum $\Sigma^{d-1}(u)$ is locally covered by a $C^{1,\alpha}$-hypersurface. 

For  $k=1,\dots,d-2$, the lower stratum $\Sigma^k(u)$ is locally covered by a $C^{1,\log^\eps}$-manifold of $k$ dimensions.

The lowest stratum $\Sigma^0(u)$ consists of isolated points. 
\end{thm} 
\section{The fully nonlinear obstacle problem}
Many ideas developed for the classical obstacle problem have been successfully modified and applied to other free boundary problems. The techniques for the regular part are especially robust, as they can be applied to problems with nonlinear operators and even integro-differential operators \cite{CSR}. 

On the other hand, the treatment for the singular set has not been as widely adapted.  The main obstruction is the reliance on monotonicity formulae. These are powerful tools, but they are also very restrictive, essentially working only for the Laplacian operator. 

To understand the regularity of the singular set in \FB problems with more general operators, it is desirable to develop tools that do no rely on monotonicity formulae.  Recently, joint with O.~Savin, we have successfully studied the singular set in the fully nonlinear obstacle problem \cite{SY}.  In the remaining pages, we highlight some ideas in this work. 
\subsection{Fully nonlinear elliptic operators}Let $\mathcal{S}_d$ denote the space of symmetric real $d$-by-$d$ matrices. A fully nonlinear operator is a function $F:\mathcal{S}_d\to\R$. For this operator to be elliptic, we further require $F$ to satisfy, for some constant $1\le\Lambda<+\infty$, the following \textit{ellipticity condition} \begin{equation}\label{Ellipticity}\frac{1}{\Lambda}\|P\|\le F(M+P)-F(M)\le\Lambda\|P\| \text{ for all $M, P\in\mathcal{S}_d$ and $P\ge 0.$ }\end{equation}

Such operators are generalizations of the Laplacian in the sense that they enjoy the maximum principle.  However, their highly nonlinear nature means that they lack any divergence structure.  In particular, no monotonicity formulae are expected for such operators. For the regularity theory of fully nonlinear elliptic operators, the standard reference is Caffarelli-Cabr\'e \cite{CC}.

Besides their importance in applications, fully nonlinear elliptic operators are very interesting objects theoretically. The absence of divergence structure pushes us to understand the deep mechanisms that make elliptic theory work. This is certainly the case in the study of the singular set in the fully nonlinear obstacle problem.
\subsection{The fully nonlinear obstacle problem} We study \textit{the fully nonlinear obstacle problem}, that is, \begin{equation}\label{FullyNonlinearObstacle}
\begin{cases}F(D^2u)=\chi_{\{u>0\}}&\\u\ge 0&
\end{cases} \text{ in $\Omega$.}
\end{equation} Here $\Omega$ is a domain in $\R^d$, and $\chi_E$ denotes the characteristic function of a set $E$. The free boundary in this problem is $\partial\PosS$.  The goal is to understand the regularity of the solution as well as the regularity of this free boundary. 

Since $F(D^2u)$ is discontinuous along the free boundary, the best possible regularity of the solution is $C^{1,1}_{loc}$.  To achieve this, we need to impose conditions on $F$ so that solutions to the equation with constant right-hand side, $F(D^2v)=1$, is better than $C^{1,1}$. By Evans-Krylov theorem, it suffices to impose the following:
$$\text{$F$ is convex.}$$It is natural to assume the obstacle itself is a solution, that is, $$F(0)=0.$$For some technical reasons, we  also assume that $F$  is $C^1$.
\subsection{Regularity of the solution and the regular part of the free boundary}The initial steps in the regularity theory of the classical obstacle problem are based  on the maximum principle. These include  the regularity of the solution, the classification of blow-up profiles, and  the regularity of the regular part of the free boundary. These can be generalized to the fully nonlinear version, since fully nonlinear operators still enjoy the maximum principle. 

This was done by Lee in his thesis \cite{L}. 

He proved the optimal regularity and non-degeneracy property of the solution, as in Theorem \ref{OptimalRegularity} and Proposition \ref{QuadraticGrowth}. The only difference is that the estimates further depend on the ellipticity constant $\Lambda$ in \eqref{Ellipticity}.

The classification of blow-up profiles is  similar to  Proposition \ref{Classification}. Due to the nonlinearity, half-space solutions are of the form $$u_0(x)=\frac{1}{2}c_e\max\{x\cdot e,0\}^2$$ for some direction $e\in\Sph$ and a constant $c_e$ depending on $e$. Parabola solutions are of the form $$u_0(x)=\frac{1}{2}x\cdot Ax$$ for a matrix $A\ge 0$ satisfying $F(A)=1.$ We define the space of parabola solutions as 
\begin{equation}\label{ParabolaSolutions}
\Qua=\{p=\frac{1}{2}x\cdot Ax:A\ge0, F(A)=1\}.
\end{equation}

The geometric characterization of regular and singular points is similar to Proposition \ref{GeometricCharacterization}. Just as in the case of the classical obstacle problem,   there is true contact in `full measure' around a regular point. Around a singular point, there is only tangential contact.

The regular part enjoys similar regularity as in Theorem \ref{CaffarelliResult}. With a fully nonlinear operator, however, the regular part is not analytic in general. It is always locally a $C^{1,\alpha}$ hypersurface. If the operator is smooth, then this can be bootstraped to $C^\infty$.

\subsection{The singular set in the \FB}The study of the regular part in the fully nonlinear obstacle problem is very similar to the case in the classical obstacle problem.  The stability of the \FB near a regular point allows arguments based entirely on the maximum principle. 

This is not the case for the singular part. As we have seen in Section \ref{SingularPartClassical}, to deal with the unstable nature of the \FB near a singular point,  all kinds of monotonicity formulae were invoked in previous studies.  These formulae are very powerful, but also restrictive, essentially working only for the Laplacian.  Since fully nonlinear operators lack any kind of divergence structure, monotonicity formulae are not expected in the fully nonlinear obstacle problem.

Consequently, for a long time, little was known about the singular part  for the fully nonlinear obstacle problem. Even the uniqueness of blow-up profile was not clear. This absence of monotonicity formulae is also an obstruction to the study of the singular set in other problems with  nonlinear operators. It is desirable to develop tools that do not rely on monotonicity formulae.  
 
 Recently, this has been achieved in a joint work with O.~Savin \cite{SY}. We proved the following:
 
 \begin{thm}\label{MainResult}
Let $u$ be a solution to the fully nonlinear obstacle problem \eqref{FullyNonlinearObstacle}. 

For $k=1,\dots,d-2$, the $k$-th stratum of the singular points, $\Sigma^k(u)$, is locally covered by a $k$-dimensional $C^{1,\log^\eps}$-manifold. 

 The top stratum, $\Sigma^{d-1}(u)$, is locally covered by a $(d-1)$-dimensional $C^{1,\alpha}$-manifold. \end{thm}  Our result essentially matches Theorem \ref{FSeResult} in Figalli-Serra \cite{FSe}, when the operator is the Laplacian.
 
 As explained in Section \ref{SingularPartClassical}, the key is to find a parabola solution $p\in\Qua$ such that  
 \begin{equation}\label{Approximation}
|u-p|\le\sigma(r)r^2 \text{ in $B_r$ for all $r>0$}
\end{equation} for some modulus of continuity $\sigma$. 

Simply knowing $\sigma(0+)=0$ implies the uniqueness of the blow-up, and gives the differentiability of the covering manifold. More details on $\sigma$ allow comparison between tangent spaces at nearby points, which leads to better regularity of the manifolds. In particular, to get $C^{1,\alpha}$-regularity, we need to show $\sigma(r)=r^\alpha$. To get $C^{1,\log^\eps}$-regularity, we need $\sigma(r)=(-\log(r))^{-\eps}.$

To get estimate \eqref{Approximation} without monotonicity formulae, we need to understand the deep mechanism that regularizes the singular part in the free boundary.  

Since the membrane only contacts the obstacle  \textit{tangentially} near a singular point,  the obstacle affects the shape of the membrane only marginally.  Intuitively,  at a singular point, the solution $u$ to the obstacle problem `almost' solves an equation with a constant right-hand side as if the obstacle was not there.  This should imply better regularity of the solution at a singular point. 

This is evident at the level of blow-up profiles. The half-space solution $\frac{1}{2}c_e\max\{x\cdot e,0\}^2$, which models the behavior  near a regular point, is $C^{1,1}$ and not better. On the other hand, any parabola solution in $\Qua$ are smooth in the entire space.  This is also clear from the estimate we need. Estimate \eqref{Approximation} implies that $u$ is $C^{2}$ at $0$, since $u$ has a second-order  expansion. Moreover, the Hessian of $u$, $D^2u$, is continuous at $0$ with $\sigma$ as its modulus of continuity.

With this observation, what we need is an improvement of regularity of the solution near a singular point. To this end, \textit{the unstable nature of the \FB is helping us, as this instability implies  strong rigidity properties of the solution.} 
\subsection{Improvement of regularity of the solution near a singular point} 
To focus on the main ideas, we assume $$F(D^2u)=\Delta u.$$ We are actually studying the classical obstacle problem. 

Even for the case of the Laplacian, this argument is interesting, as it explains the regularizing mechanism of the singular set. So far, this mechanism has been hidden behind monotonicity formulae. In particular, our argument is entirely based on the maximum principle. 

The desired estimate \eqref{Approximation} says that we can find a parabola solution $p$, which approximates our solution $u$ better and better at smaller and smaller scales. How this approximation improves as $r\to0$ is dictated by the modulus of continuity $\sigma.$

Without monotonicity formulae, it is difficult to get such approximation at all small scales at once. Instead, we discretize and achieve this estimate inductively.  One step in this discretized version is the following:

\begin{lem}[One step improvement]\label{OneStep}Suppose that $u$ solves the fully nonlinear obstacle problem \eqref{FullyNonlinearObstacle} with $0\in\Singu.$

There are constants $\eps_0, r_0\in (0,1/2)$ such that the following is true:

 If for a parabola solution $p\in\Qua$ we have $$|u-p|<\eps \text{ in $B_1$ for some $\eps<\eps_0$,}$$ then there is  $p'\in\Qua$ such that $$|u-p'|<\eps'r_0^2 \text{ in $B_{r_0}$ for some $\eps'<\eps.$}$$
\end{lem} 

That is, if our solution is well-approximated by a parabola solution $p$ at the unit scale,  then we can find a (possibly different) parabola solution $p'$ that improves this approximation at a smaller scale $r_0.$

Once we establish this lemma, we look at the rescaled solution $u_{1}(x)=\frac{1}{r_0^2}u(r_0x)$.  It satisfies $$|u_1-p'|\le \eps'\text{ in $B_1$.}$$ Consequently, we can apply the lemma to $u_1$, with $p$ and $\eps$ replaced by $p'$ and $\eps'$. We  iterate,  producing a sequence of parabola solutions $p_k$ and a sequence of error bounds $\eps_k.$ Each iteration improves the approximation by $\eps-\eps'.$ If this improvement is significant, then the sequence $p_k$ converges to a limiting parabola $p_\infty,$ which gives the desired approximation \eqref{Approximation}. The modulus of continuity $\sigma$ is dictated by the decay $\eps\to\eps'.$

To get this one step improvement, we need to consider two cases, depending on the initial parabola solution $p$, in particular, on the second largest eigenvalue of $D^2p.$

Up to a rotation, we assume that $p$ is of the form \begin{equation}\label{StandardParabola}p(x)=\frac{1}{2}\sum a_jx_j^2\end{equation} with $a_1\ge a_2\ge\dots\ge a_d\ge 0.$

\subsubsection{Case 1: The initial parabola has only one large eigenvalue}To be concrete, let's assume $$p(x)=\frac{1}{2}x_1^2.$$ Since this parabola is strictly convex and strictly positive away from the hyperplane $\{x_1=0\}$, we  expect the contact set $\{u=0\}$ to concentrate along this hyperplane. 

Indeed,  from the comparison $u\ge\frac{1}{2}x_1^2-\eps$, it follows that the contact set is contained in the strip $\{|x_1|<C\eps^{1/2}\}$. By a barrier, we can further localize the contact set in the strip $\{|x_1|<C\eps\}.$ Note that along the boundary of this strip the parabola satisfies $|\nabla p|\le C\eps.$

Define the normalized solution $$\hat{u}_\eps=\frac{1}{\eps}(u-p).$$ This function is harmonic outside a strip of width $\eps$. 
The uniform Lipschitz bound on $p/\eps$ gives enough compactness of the family $\{\hat{u}_\eps\}$ when $\eps\to0.$ Up to a subsequence, we have $$\hat{u}_\eps\to\hat{u}_0$$ for some function $\hat{u}_0$. 

If we can show that \textit{this limit $\hat{u}_0$ is $C^2$ at the origin}, then  its second order expansion at $0$, $q$, satisfies  the following 
\begin{align*}
u&=p+\eps\hat{u}_\eps\\&=p+\eps\hat{u}_0+\eps o(1) \text{ as $\eps\to0$}\\&=p+\eps q+\eps o(r_0^2)+\eps o(1) \text{ in $B_{r_0}$  as $\eps\to0$.}
\end{align*}
We still need some modifications, but the parabola $p'=p+\eps q$ satisfies $$|u-p'|<\frac{1}{2}\eps \text{ in $B_{r_0}$}$$ if $r_0$ and $\eps$ are small. 

We get the one step improvement as in Lemma \ref{OneStep} with $\eps'=\frac{1}{2}\eps.$

The challenge is to show that the limit $\hat{u}_0$ is $C^2$ at the origin.  Since each $\hat{u}_\eps$ solves the obstacle problem  with an obstacle concentrated in a strip $\{|x_1|<\eps\}$,  it can be shown that the limit  $\hat{u}_0$ solves \textit{the thin obstacle problem} \cite{AC}. In general, solutions to this problem is only Lipschitz. However, with the assumption $0\in\Singu$, we can use the instability of the \FB to remove the low homogeneities of $\hat{u}_0$. This gives the desired $C^2$-regularity. 

This gives Lemma \ref{OneStep} for the particular parabola $p=\frac{1}{2}x_1^2.$ A similar argument works for other parabola solutions with only one large eigenvalue.  We have 

\begin{lem}\label{CaseOne}\label{QuadraticApproxTopStratum}
Suppose for some given $\kappa>0$, we have, for some $p\in\Qua$ with $a_2\le\kappa\eps,$ $$|u-p|<\eps \text{ in $B_1$}$$and $$0\in\Singu.$$

There are constants $\eps_0,r_0\in(0,1/2)$, depending on $\kappa$, such that if $\eps<\eps_0$, then $$|u-p'|<\frac{1}{2}\eps \text{ in $B_{r_0}$}$$ for some $p'\in\Qua$.\end{lem}
Here $a_2$ is the second largest eigenvalue of $D^2p$ as in  \eqref{StandardParabola}.

If in the iterations we always end up in the case described in this lemma, then we have a sequence of parabola solutions $p_k$ satisfying $$|D^2p_k-D^2p_{k+1}|\le C(1/2)^k.$$ This geometric decay implies the convergence $p_k\to p_\infty.$ The limit satisfies $$|u-p_\infty|<Cr^{2+\alpha} \text{ in $B_r$.}$$  This is consistent with the $C^{1,\alpha}$-regularity of the covering manifold for the top stratum.

\subsubsection{Case 2: The initial parabola has more than one large eigenvalues}However, there is no guarantee that the updated parabola $p'$ in Lemma \ref{CaseOne} still satisfies $a_2\le\kappa\eps$. Consequently, we need a similar improvement when the initial parabola solution $p$ has more than one large eigenvalues, that is, when $p$ is of the form $$p=\frac{1}{2}\sum a_jx_j^2$$ with $a_1\ge a_2\ge\kappa\eps.$

In this case, $p$ is strictly positive away from the set $\{x_1=x_2=0\}$. We expect the contact set $\{u=0\}$ to be concentrated around $\{x_1=x_2=0\}$. This set is of higher co-dimensions.  In particular, there are no barrier functions. 

On the other hand, since the contact set is so `thin', it is more natural to use the solution to \textit{the unconstraint problem} $$\begin{cases}\Delta h=1&\text{ in $B_1$}\\h=u&\text{ along $\partial B_1$}\end{cases}$$to approximate our solution. That is, we want to take the second order expansion of $h$ at the origin as the next parabola solution $p'$. However, without the  constraint $h\ge0$,  this second order expansion is  not a parabola solution in $\Qua.$ 

Another subtlety is that there is no geometric decay $\eps\to\frac{1}{2}\eps$ as in the previous case.  Here we are considering the case when the contact set is of higher co-dimensions, for which the covering manifold is of $C^{1,\log^\eps}$-regularity.  The corresponding decay  is of the form $\eps\to\eps-\eps^M$ for some large  power $M$. This  slow decay is not enough to control the distance between consecutive parabolas appearing in the iterations.

These delicate issues require very detailed analysis in this case. We do not have enough space to produce the arguments here. The interested reader should consult Savin-Yu \cite{SY}.


\end{document}